\def\PaperDate{2004/05/17}
\documentclass[12pt]{article}
  \usepackage[latin1]{inputenc}
  \usepackage[american]{babel}
  \usepackage[USletter]{buxlayout}
  \usepackage[GlobalNumbered]{buxmath}
  \usepackage{datum}
  \usepackage{xypic}
  \usepackage[dvips]{graphicx}
  \usepackage[noquotes]{bux2ref}
  \InputIfFileExists{basic_notation.tex}{}{}

  \usepackage{ulem}
  
  \newcommand{\Pair}[2]{(#1,#2)}
  \newcommand{\crossprod}{\times}

  \newvariable{\TheTrivialElement}{\Number{1}}
  \newvariable{\TheLastIndex}{r}
  \newvariable{\TheVariable}{x}
  \newvariable{\TheWord}{w}
  \newvariable{\AltWord}{v}
  \newvariable{\ThirdWord}{u}
  \newvariable{\GriOne}{\alpha}
  \newvariable{\GriTwo}{\beta}
  \newvariable{\GriThree}{\gamma}
  \newvariable{\TheGroup}{G}
  \newcommand{\subgroup}{\leq}
  \newvariable{\TheSubgroup}{H}

  \newvariable{\AltGroupElement}{h}
  \newvariable{\AltGroup}{H}
  \newcommand{\Choose}[2]{\left(\begin{array}{@{}c@{}}#1\\#2\end{array}\right)}
  \newvariable{\DihedralGroup}{D}
  \newcommand{\FloorOf}[1]{\left\lfloor #1 \right\rfloor}
  \newvariable{\ImUnit}{i}
  \newvariable{\TheBall}{B}
  \newvariable{\TheConstant}{C}
  \newvariable{\TheEpsilon}{\varepsilon}
  \newvariable{\TheGeneratingSet}{\Sigma}
  \newvariable{\TheGenerator}{x}
  \newvariable{\TheGroupElement}{g}
  \newvariable{\TheGrowthRate}{\lambda}
  \newvariable{\TheHomomorphism}{\phi}
  \newvariable{\TheLength}{\ell}
  \newvariable{\TheNumber}{n}
  \newcommand{\ThePositiveReals}{\RRR^{+}}
  \newvariable{\TheRadius}{r}
  \newvariable{\TheShiftConst}{K}
  \newvariable{\TheSign}{\varepsilon}

  \newvariable{\AltVariable}{y}
  \newvariable{\TheComplex}{z}
  \newcommand{\FundamentalGroupOf}[2][\sadkjal]{\pi_1(#2,#1)}
  \newvariable{\ThePoint}{x}
  \newvariable{\TheLevel}{s}
  \newvariable{\TheAutomorphism}{\alpha}
  \newcommand{\CardOf}{\#}
  \newvariable{\RiemannSphere}{\overline{\CCC}}
  \newcommand{\ImgOf}[1]{\operatorname{IMG}\left(#1\right)}
  \newvariable{\TheBound}{M}
  \newvariable{\BadStrings}{V}
  \newcommand{\pref}[1]{(\ref{#1})}
  \newcommand{\Slot}{\textvisiblespace}  
    
  \newvariable{\TheOrbit}{\Oka}
  \newvariable{\TheConstTerm}{c}
  
  \newvariable{\TheBasePoint}{*}
  \newvariable{\TheLoop}{\gamma}
  \newvariable{\TheLift}{\tilde{\TheLoop}}
  \renewcommand{\Left}{\mathsf{L}}
  \renewcommand{\Right}{\mathsf{R}}
  \global\let\TheRootedTree\undefined
  \newvariable{\TheRootedTree}{T^{(2)}}
  \newvariable{\TheSwap}{\sigma}
  \newvariable{\MaybeSwap}{\varepsilon}
  \newvariable{\ThePercentage}{\eta}
  \newvariable{\AltPercentage}{\kappa}
  \newvariable{\TheGrowth}{\operatorname{gr}}
  \newvariable{\GroupG}{\mathcal{G}}
  \newvariable{\GroupH}{\mathcal{H}}
  \newvariable{\GroupI}{\mathcal{I}}
  \newvariable{\GroupJ}{\mathcal{J}}
  \newvariable{\TheAutom}{\alpha}
  \newvariable{\AltAutom}{\beta}
  \newvariable{\ThrAutom}{\gamma}
  \newvariable{\TheVar}{x}
  \newvariable{\AltVar}{y}
  \newvariable{\ThrVar}{z}
  \newvariable{\eee}{\mathrm{e}}
  \global\let\CyclicGroup\undefined
  
  \newvariable{\CyclicGroup}{\ZZZ}
  \newvariable{\TheBadWords}{b}
  \newvariable{\TheBadElements}{B}
  \renewcommand{\ExpOf}[1]{\eee^{#1}}
  \newvariable{\ThePostcriticalSet}{\mathcal{P}}
  \newvariable{\ThePolynomial}{f}
  \newvariable{\AltPolynomial}{P}
  \newvariable{\CoverLength}{\tilde{\TheLength}}
  \newvariable{\TheLanguage}{\TheGeneratingSet^*}
  \newvariable{\AltRadius}{\TheRadius'}
  \newvariable{\TheProportion}{p}
  \newcommand{\TheNonnegativeReals}{\RRR^+_0}
  \newcommand{\AutOf}[1]{\Aut(#1)}
  \newcommand{\ClosureOf}[1]{\overline{#1}}
  \newvariable{\TheIndex}{i}
  \newvariable{\AltIndex}{j}
\begin{document}
  \title{On the Growth of Iterated Monodromy Groups\footnotetext{Department of Mathematics, Cornell University, Ithaca, NY 14850.}\footnotetext{2000 AMS Classification: {\bf 20F65, 37F20.}}}
  \author{Kai-Uwe Bux \and Rodrigo P\'erez\footnote{Research supported by NSF Postdoctoral Fellowship, grant DMS-0202519.}}
  \date{\datum{\PaperDate}}

  \maketitle

  \begin{abstract}
    Nekrashevych conjectured that the iterated monodromy groups of quadratic
    polynomials with preperiodic critical orbit have intermediate growth. We
    illustrate some of the difficulties that arise in attacking this
    conjecture and prove subexponential growth for the iterated monodromy
    group of $z^2+i$. This is the first non-trivial example supporting the
    conjecture.
  \end{abstract}

  \noindent
    The iteration of a quadratic polynomial
    \(
      \ThePolynomial =
      \ThePolynomialOf[\TheConstTerm]{\TheComplex}
      =
      \ToThePower{\TheComplex}{\Two} + \TheConstTerm
    \)
    describes a dynamical system in $\RiemannSphere$.
    The behavior of this system is ruled
    by the geometry of the orbit
    \[
      \TheOrbit =
      \TheOrbit[\TheConstTerm] :=
      \SetOf{
        \ThePolynomialOf{\Zero}
        ,
        \ThePolynomialOf{\ThePolynomialOf{\Zero}}
        ,\ldots,
        \ThePolynomialOf[][\TheIndex]{\Zero}
        ,\ldots
      }
    \]
    of its unique critical point \cite{Milnor:1999}.
    V.~Nekrashevych \cite{Nekrashevych:2003} associates to each
    such system a group of automorphisms of the infinite binary rooted
    tree $\TheRootedTree$. In Section~\ref{SEC::IMG_Construction},
    we will sketch the construction
    of this group known as the
    \notion{iterated monodromy group} of $\ThePolynomial$ denoted
    by $\ImgOf{\ThePolynomial}$.
    This note addresses the following conjecture:
    \begin{conj}[Nekrashevych]\label{Conjecture}
      Suppose the critical orbit of $\ThePolynomial$ is postcritically
      finite, i.e., the orbit of the critical point
      $\Zero$ is finite and does not contain $\Zero$.
      Then $\ImgOf{\ThePolynomial}$
      has intermediate growth.
    \end{conj}
    We want to illustrate some of the difficulties that arise in
    attacking this conjecture. Our plan is to present three examples
    $\GroupG$, $\GroupH$ and $\GroupI$ of finitely generated subgroups
    of $\AutOf{\TheRootedTree}$, all of which have subexponential
    growth. The group $\GroupG$ is the First Grigorchuk group; it was
    the first known example of a group of intermediate growth
    \cite{Grigorchuk:1983}. The group $\GroupH$ belongs to the family
    of groups of intermediate growth studied in
    \cite{Grigorchuk:1984}. Our proofs of subexponential growth for
    $\GroupG$ and $\GroupH$ are designed to illustrate the use of
    Proposition \ref{BasicTool}. We apply these ideas to prove
    subexponential growth on $\GroupI := \ImgOf{\TheComplex[][\Two] +
    \ImUnit}$, thus providing the first non-trivial example to support
    Conjecture \ref{Conjecture}.

    We want to express our gratitude to R. Grigorchuk, J. Hubbard and J.
    Smillie for their interest in this project. In particular, we thank R.
    Grigorchuk for his comments on a previous version and for supplying a
    reference.

  \section{The Iterated Monodromy Group}\label{SEC::IMG_Construction}
    The \notion{postcritical set} of $\ThePolynomial$ is
    the closure
    \(
      \ClosureOf{\TheOrbit}
    \)
    of the critical orbit
    \(
      \TheOrbit=
      \SetOf[{
        \ThePolynomialOf[][\TheIndex]{\Zero}
      }]{
        \TheIndex \geq \One
      }
      .
    \)
    Note that
    \(
      \ThePolynomialOf{\ClosureOf{\TheOrbit}}
      \subseteq
      \ClosureOf{\TheOrbit}
      ,
    \)
    whence
    \(
      \ThePolynomialOf[][-\One]{
        \RiemannSphere\setminus\ClosureOf{\TheOrbit}
      }
      \subseteq
      \RiemannSphere\setminus\ClosureOf{\TheOrbit}
    \)
    and
    \(
      \ThePolynomial
      \mapcolon 
      \ThePolynomialOf[][-\One]{
        \RiemannSphere\setminus\ClosureOf{\TheOrbit}
      }
      \rightarrow
      \RiemannSphere\setminus\ClosureOf{\TheOrbit}
    \)
    is a $\Two$ to $\One$ covering map.

    \begin{figure}
      \begin{center}
      \setlength{\unitlength}{0.7cm}
      \begin{picture}(18,11)(-2,-1)
        \put(-0.2,-1.5){\includegraphics{img.1}}
        \put(7,0){\lap[10pt,0pt]{$\TheBasePoint[\emptyset]$}}
        \put(3,3){\lap[10pt,0pt]{$\TheBasePoint[\Left]$}}
        \put(11,3){\lap[10pt,0pt]{$\TheBasePoint[\Right]$}}
        \put(1,6){\lap[10pt,0pt]{$\TheBasePoint[\Left\Left]$}}
        \put(5,6){\lap[10pt,0pt]{$\TheBasePoint[\Left\Right]$}}
        \put(9,6){\lap[10pt,0pt]{$\TheBasePoint[\Right\Left]$}}
        \put(13,6){\lap[10pt,0pt]{$\TheBasePoint[\Right\Right]$}}
      \end{picture}
      \end{center}
      \caption{The monodromy action on $\TheRootedTree$%
      \label{FIG::IMA}}
    \end{figure}
    
    Fix a base-point $\TheBasePoint =
    \TheBasePoint[\emptyset]\in\RiemannSphere\setminus\ClosureOf{\TheOrbit}$.
    This point has two $\ThePolynomial$-preimages
    $\TheBasePoint[\Left]$ and $\TheBasePoint[\Right]$.
    For any loop $\TheLoop$ in
    $\RiemannSphere\setminus\ClosureOf{\TheOrbit}$ based
    at $\TheBasePoint$, the homotopy class
    \(
      \EquivClassOf{\TheLoop}\in
      \FundamentalGroupOf[\TheBasePoint]{\RiemannSphere\setminus\ClosureOf{\TheOrbit}}
    \)
    determines whether $\TheLoop$
    lifts to two loops
    (based at $\TheBasePoint[\Left]$ and $\TheBasePoint[\Right]$)
    or to two paths, one from $\TheBasePoint[\Left]$ to
    $\TheBasePoint[\Right]$ and one in the opposite direction.
    The \notion{monodromy action} of
    \(
      \EquivClassOf{\TheLoop}\in
      \FundamentalGroupOf[\TheBasePoint]{\RiemannSphere\setminus\ClosureOf{\TheOrbit}}
    \)
    on the set
    $\SetOf{\TheBasePoint[\Left],\TheBasePoint[\Right]}$ is the map
    that takes
    \(
      \ThePoint\in\SetOf{\TheBasePoint[\Left],\TheBasePoint[\Right]}
    \)
    to the endpoint of the lift of $\TheLoop$ based at $\ThePoint$.
    
    Note that both $\TheBasePoint[\Left]$ and $\TheBasePoint[\Right]$ have
    two $\ThePolynomial$-preimages:
    $\TheBasePoint[\Left\Left]$, $\TheBasePoint[\Left\Right]$ and
    $\TheBasePoint[\Right\Left]$, $\TheBasePoint[\Right\Right]$ respectively.
    Since
    \(
      \ThePolynomial \mapcolon
      \ThePolynomialOf[][-\One]{\RiemannSphere\setminus\ClosureOf{\TheOrbit}}
      \rightarrow
      \RiemannSphere\setminus\ClosureOf{\TheOrbit}
    \)
    is a two sheeted covering map, all these four points are different.
    The loop $\TheLoop$ has four different $\ThePolynomial[][\Two]$-lifts,
    one starting at each of the points
    \(
      \TheBasePoint[\Left\Left],
      \TheBasePoint[\Left\Right],
      \TheBasePoint[\Right\Left],
      \TheBasePoint[\Right\Left]
      .
    \)
    Sending each of these four points to the endpoint of the corresponding
    $\ThePolynomial[][\Two]$-lift of $\TheLoop$ defines an action of
    $\FundamentalGroupOf[\TheBasePoint]{\RiemannSphere\setminus\ClosureOf{\TheOrbit}}$ on the set
    \(
      \ThePolynomialOf[][-\Two]{\TheBasePoint} =
      \SetOf{
        \TheBasePoint[\Left\Left],
        \TheBasePoint[\Left\Right],
        \TheBasePoint[\Right\Left],
        \TheBasePoint[\Right\Left]
      }
      .
    \)
    Iterating this procedure, we define the
    \notion{monodromy action} of
    $\FundamentalGroupOf[\TheBasePoint]{\RiemannSphere\setminus\ClosureOf{\TheOrbit}}$
    on the sets
    \(
      \ThePolynomialOf[][-\TheLevel]{\TheBasePoint}
      .
    \)
    Figure~\ref{FIG::IMA} illustrates these actions.
    There, we have arranged the sets
    \(
      \ThePolynomialOf[][-\TheLevel]{\TheBasePoint}
    \)
    into an infinite binary rooted tree. The edges from the vertices
    at level $\TheLevel$ to their level $\TheLevel-\One$ parents correspond
    to the map
    \(
      \ThePolynomial \mapcolon
      \ThePolynomialOf[][-\TheLevel]{\TheBasePoint}
      \rightarrow
      \ThePolynomialOf[][-\Parentheses{\TheLevel-\One}]{\TheBasePoint}
      .
    \)
    \begin{rem}
      Figure~\ref{FIG::IMA} is slightly misleading:
      Although we have drawn various
      lifts of the loop $\TheLoop$, the abstract tree does not lie inside
      $\RiemannSphere\setminus\ClosureOf{\TheOrbit}$.
      Indeed, vertices of different levels might correspond
      to the same point in $\RiemannSphere\setminus\ClosureOf{\TheOrbit}$.
      This occurs exactly when
      $\TheBasePoint$ is a periodic point of $\ThePolynomial$.

      To avoid this ambiguity, we transfer the monodromy action to an
      abstract binary tree, still denoted by $\TheRootedTree$, with vertices
      labelled by words in
      \(
        \SetOf{\Left,\Right}
      \)
      in accordance to the indices used above for the preimages of
      $\TheBasePoint$.
    \end{rem}
    Note that the monodromy actions of
    \(
      \FundamentalGroupOf[\TheBasePoint]{
        \RiemannSphere\setminus\ClosureOf{\TheOrbit}
      }
    \)
    on the different levels are compatible with the edges in the
    tree $\TheRootedTree$. The \notion{iterated monodromy action}
    is the action of 
    \(
      \FundamentalGroupOf[\TheBasePoint]{
        \RiemannSphere\setminus\ClosureOf{\TheOrbit}
      }
    \)
    on $\TheRootedTree$ thus induced.
    The \notion{iterated monodromy group} of $\ThePolynomial$ is the image
    of the iterated monodromy action in $\AutOf{\TheRootedTree}$. It
    will be denoted by $\ImgOf{\ThePolynomial}$.
    \begin{observation}
      When indexing the iterated preimages of $\TheBasePoint$,
      we make choices at each step about which preimage we
      call $\Left$ and which
      one we call $\Right$. Any two labellings of the iterated preimages
      are conjugate via an automorphism of $\TheRootedTree$. Thus the
      iterated monodromy group is well defined up to conjugation in
      $\AutOf{\TheRootedTree}$.
    \end{observation}

  \section{Automorphisms of the Infinite Binary
           Tree \boldmath$\TheRootedTree$}
    The infinite binary rooted tree $\TheRootedTree$ has two
    subtrees connecting to the root vertex. We call
    them $\TheRootedTree[\Left]$ (left) and
    $\TheRootedTree[\Right]$ (right). Both subtrees are binary
    infinite rooted trees in their own right. Fix isomorphisms
    \(
      \TheRootedTree[\Left] \cong \TheRootedTree
    \)
    and
    \(
      \TheRootedTree[\Right] \cong \TheRootedTree
      .
    \)
    Using these identifications, any two automorphisms
    $\TheAutomorphism[\Left],\TheAutomorphism[\Right]
    \in\AutOf{\TheRootedTree}$
    can be combined to define an automorphism
    $\TheAutomorphism :=
    \Pair{\TheAutomorphism[\Left]}{\TheAutomorphism[\Right]}
    \in\AutOf{\TheRootedTree}$: the automorphism $\TheAutomorphism$
    preserves each of the two subtrees and
    acts in $\TheRootedTree[\Left]$ as $\TheAutomorphism[\Left]$ and
    in $\TheRootedTree[\Right]$ as $\TheAutomorphism[\Right]$.

    Moreover, the isomorphism
    \(
      \TheRootedTree[\Left] \cong \TheRootedTree
      \cong
      \TheRootedTree[\Right]
    \)
    provides a distinguished involution $\TheSwap\in\AutOf{\TheRootedTree}$,
    called the \notion{swap},
    that interchanges $\TheRootedTree[\Left]$ and $\TheRootedTree[\Right]$.
    Note that the swap interacts nicely with the \notion{pair notation}
    from above:
    \[
      \TheSwap
      \Pair{\TheAutomorphism[\Left]}{\TheAutomorphism[\Right]}
      \TheSwap[][-\One]
      =
      \Pair{\TheAutomorphism[\Right]}{\TheAutomorphism[\Left]}
      .
    \]
    It is easy to see that this defines a wreath product decomposition
    \[
      \AutOf{\TheRootedTree}
      =
      \Parentheses{
        \AutOf{\TheRootedTree}
        \crossprod
        \AutOf{\TheRootedTree}
      }
      \lprod
      \GroupPresented{\TheSwap}
    \]
    which allows us to represent automorphisms
    of $\TheRootedTree$ pictorially:
    \[
      \Pair{\TheAutomorphism[\Left]}{\TheAutomorphism[\Right]}
      \TheSwap
      =
      \begin{array}{c}
        \xymatrix@-15pt{
          {\TheAutomorphism[\Left]} & & {\TheAutomorphism[\Right]} \\
                 & \TheSwap \ar@{-}[ul] \ar@{-}[ur]  & \\
        }
      \end{array}
      \qquad\quad
      \Pair{\TheAutomorphism[\Left]}{\TheAutomorphism[\Right]}
      =
      \begin{array}{c}
        \xymatrix@-15pt{
          {\TheAutomorphism[\Left]} & & {\TheAutomorphism[\Right]} \\
                 & \TheTrivialElement \ar@{-}[ul] \ar@{-}[ur]  & \\
        }
      \end{array}
    \]
    \begin{lemma}
      For $\TheIndex=\One,\ldots,\TheLastIndex$,
      let $\TheVariable[\TheIndex]$ be variables with values in
      $\AutOf{\TheRootedTree}$, let $\AltWord[\TheIndex]$ and
      $\ThirdWord[\TheIndex]$ be words in $\TheSwap$ and the variables.
      For each index $\TheIndex$, fix
      $\MaybeSwap[\TheIndex]\in\SetOf{\TheTrivialElement,\TheSwap}$.
      Then the system
      \begin{eqnarray*}
        \TheVariable[\One] &=& \Pair{\AltWord[\One]}{\ThirdWord[\One]} \MaybeSwap[\One] \\
                           &\vdots&  \\
        \TheVariable[\TheLastIndex] &=& \Pair{\AltWord[\TheLastIndex]}{\ThirdWord[\TheLastIndex]} \MaybeSwap[\TheLastIndex]
      \end{eqnarray*}
      has a unique solution.
    \end{lemma}
    Instead of a formal proof, we give a convincing example:
    \begin{example}
      Consider the system
      \begin{eqnarray*}
        \GriOne &=& \Pair{\TheSwap}{\GriTwo} \\
        \GriTwo &=& \Pair{\TheSwap}{\GriThree} \\
        \GriThree &=& \Pair{\TheTrivialElement}{\GriOne}
      \end{eqnarray*}
      First, we obtain information about how these three automorphisms
      act on the first layer in $\TheRootedTree$:
      \[
        \TheAutom
        =
        \begin{array}{c}
          \xymatrix@-15pt{
            \TheSwap & & \AltAutom \\
                   & \TheTrivialElement \ar@{-}[ul] \ar@{-}[ur]  & \\
          }
        \end{array}
        \quad
        \AltAutom
        =
        \begin{array}{c}
          \xymatrix@-15pt{
            \TheSwap & & \ThrAutom \\
                   & \TheTrivialElement \ar@{-}[ul] \ar@{-}[ur]  & \\
          }
        \end{array}
        \quad
        \ThrAutom
        =
        \begin{array}{c}
          \xymatrix@-15pt{
            \TheTrivialElement & & \TheAutom \\
                   & \TheTrivialElement \ar@{-}[ul] \ar@{-}[ur]  & \\
          }
        \end{array}
      \]
      Now, we substitute and extend the picture to
      the second layer:
      \[
        \TheAutom
        =
        \begin{array}{c}
          \xymatrix@-15pt{
                     & \TheSwap & & \ThrAutom \\
            \TheSwap & & \TheTrivialElement\ar@{-}[ul] \ar@{-}[ur] \\
                   & \TheTrivialElement \ar@{-}[ul] \ar@{-}[ur]  & & \\
          }
        \end{array}
        \quad
        \AltAutom
        =
        \begin{array}{c}
          \xymatrix@-15pt{
                     & \TheTrivialElement & & \TheAutom \\
            \TheSwap & & \TheTrivialElement\ar@{-}[ul] \ar@{-}[ur] \\
                   & \TheTrivialElement \ar@{-}[ul] \ar@{-}[ur]  & & \\
          }
        \end{array}
        \quad
        \ThrAutom
        =
        \begin{array}{c}
          \xymatrix@-15pt{
                     & \TheSwap & & \AltAutom \\
            \TheTrivialElement & & \TheTrivialElement\ar@{-}[ul] \ar@{-}[ur] \\
                   & \TheTrivialElement \ar@{-}[ul] \ar@{-}[ur]  & & \\
          }
        \end{array}
      \]
      We continue and see that the system of equations determines
      $\GriOne$, $\GriTwo$, and $\GriThree$ completely.
    \end{example}
    Finally, we introduce a bit of notation.
    For any subgroup $\TheGroup\subgroup\AutOf{\TheRootedTree}$, we
    let $\TheGroup[\TheLevel]$ denote the subgroup of those
    elements in $\TheGroup$ that fix pointwise the set of vertices
    within distance $\TheLevel$ from the root. An automorphism can be
    written as a pair
    $\Pair{\TheAutomorphism[\Left]}{\TheAutomorphism[\Right]}$
    if and only if it is in $\TheGroup[\One]$. We call those elements
    of $\TheGroup$ \notion{even}. Note that
    $\TheGroup[\TheLevel]$ is a normal subgroup of finite index
    in $\TheGroup$. The even elements from a subgroup
    of index at most $\Two$.

  \section{Group Growth via Weights}
    \begin{Def}
      Consider a group $\TheGroup$ with a fixed finite generating set
      $\TheGeneratingSet$. Any map
      \[
        \CoverLength \mapcolon
        \TheGeneratingSet
        \rightarrow
        \ThePositiveReals,
      \]
      assigning a strictly positive weight to each generator, extends
      to a \notion{length function} on the set $\TheLanguage$ of words
      over the alphabet $\TheGeneratingSet\union\TheGeneratingSet[][-\One]$:
      \begin{eqnarray*}
        \CoverLength \mapcolon \TheLanguage & \rightarrow & \TheNonnegativeReals \\
        \TheWord =
          \TheGenerator[\One][{\TheSign[\One]}]
          \TheGenerator[\Two][{\TheSign[\Two]}]
          \cdots
          \TheGenerator[\TheLastIndex][{\TheSign[\TheLastIndex]}]
        & \mapsto &
          \Sum[\TheIndex=\One][\TheLastIndex]{
            \CoverLengthOf{\TheGenerator[\TheIndex]}
          }
          \qquad
          \text{(here\ }
          \TheSign[\TheIndex]=\pm\One\text{)}
      \end{eqnarray*}
      This length descends to a length function
      $\TheLength$ on $\TheGroup$ as follows:
      \begin{eqnarray*}
        \TheLength \mapcolon \TheGroup & \rightarrow & \TheNonnegativeReals \\
        \TheGroupElement & \mapsto &
          \Min{
            \SetOf[{ \CoverLengthOf{\TheWord} }]{
            \TheWord \text{\ represents\ } \TheGroupElement
          }}
      \end{eqnarray*}
    \end{Def}
    \begin{observation}
      The map $\TheLength \mapcolon \TheGroup \rightarrow \TheNonnegativeReals$
      satisfies the following conditions:
      \begin{enumerate}
        \item
          For any group element $\TheGroupElement\in\TheGroup$, we have
          \(
            \TheLengthOf{\TheGroupElement} \geq \Zero
            .
          \)
        \item
          We have $\TheLengthOf{\TheGroupElement}=\Zero$ if and only
          if $\TheGroupElement=\TheTrivialElement$.
        \item
          For any two group elements $\TheGroupElement,\AltGroupElement
          \in \TheGroup$, the inequality
          \(
            \TheLengthOf{\TheGroupElement\AltGroupElement}
            \leq
            \TheLengthOf{\TheGroupElement} + \TheLengthOf{\AltGroupElement}
          \)
          holds.
        \item
          For any radius $\TheRadius$, the set
          \(
            \TheBallOf[\TheLength]{\TheRadius}
            :=
            \SetOf[
              \TheGroupElement\in\TheGroup
            ]{
              \TheLengthOf{\TheGroupElement} \leq \TheRadius
            }
          \)
          is finite. We call this set
          the \notion{ball of radius $\TheRadius$}.
      \end{enumerate}
    \end{observation}
    \begin{Def}
      The \notion{growth
      function} $\TheGrowth[\CoverLength]$ associated to $\TheGeneratingSet$
      and
      $\CoverLength\mapcolon\TheGeneratingSet\rightarrow\ThePositiveReals$
      is defined as the ``combinatorial volume'' of the ball of radius
      $\TheRadius$:
      \[
        \TheGrowthOf{\TheRadius} = 
        \TheGrowthOf[\CoverLength]{\TheRadius}
        :=
        \CardOf{
            \SetOf[
              \TheGroupElement\in\TheGroup
            ]{
              \TheLengthOf{\TheGroupElement} \leq \TheRadius
            }
        }
        .
      \]
      The growth function is almost submultiplicative: There is a
      constant $\TheConstant$ such that
      \[
        \TheGrowthOf{\TheRadius+\AltRadius}
        \leq
        \TheConstant
        \TheGrowthOf{\TheRadius}
        \TheGrowthOf{\AltRadius}
        .
      \]
      It follows that the \notion{growth rate}
      \[
        \TheGrowthRate = \Lim[\TheRadius\rightarrow\infty]{
          \ToThePower{\Parentheses{
            \TheGrowthOf{\TheRadius}
          }
          }{1/\TheRadius}
        }
      \]
      exists. We say that $\TheGroup$ has \notion{exponential growth}
      if $\TheGrowthRate>\One$ and \notion{subexponential growth}
      if $\TheGrowthRate\leq\One$.
    \end{Def}
    \begin{rem}
      The growth rate depends on the choices of $\TheGeneratingSet$ and
      $\CoverLength\mapcolon\TheGeneratingSet\rightarrow\ThePositiveReals$.
      However, any other choice will yield a Lipschitz equivalent length
      function on $\TheGroup$. It follows that whether $\TheGroup$ has
      exponential or subexponential growth does not depend on these choices.
    \end{rem}
    The following proposition is modeled
    upon the standard proof of subexponential growth for the
    First Grigorchuk group.
    \begin{prop}\label{BasicTool}
      Let $\AltGroup$ be a finite index subgroup of the finitely
      generated group $\TheGroup$, and let $\TheLength$ be a
      length function on $\TheGroup$ as above.
      Suppose there exists numbers
      $\ThePercentage\in\RightOpenInterval{\Zero}{\One}$,
      $\TheProportion\in\LeftOpenInterval{\Zero}{\One}$,
      $\TheShiftConst\geq\Zero$, and an injective
      homomorphism
      \begin{eqnarray*}
        \TheHomomorphism \mapcolon
        \AltGroup
        & \rightarrow &
        \overbrace{%
          \TheGroup \crossprod \cdots \crossprod \TheGroup
        }^{\TheNumber \text{\ factors}}
        \\
        \AltGroupElement
        & \mapsto &
        \TupelOf{
          \TheHomomorphismOf[\One]{\AltGroupElement}
          ,\ldots,
          \TheHomomorphismOf[\TheNumber]{\AltGroupElement}
        }
      \end{eqnarray*}
      such that the following condition holds:
      \begin{quote}
        For each $\TheRadius$, the proportion of all elements
        in
        \(
          \SetOf[
            \AltGroupElement\in\AltGroup
          ]{
            \TheLengthOf{\AltGroupElement} \leq \TheRadius
          }
        \)
        satisfying
        \[
          \Sum[\TheIndex=\One][\TheNumber]{
            \TheLengthOf{
              \TheHomomorphismOf[\TheIndex]{\AltGroupElement}
            }
          }
          \leq
          \ThePercentage \TheRadius + \TheShiftConst
        \]
        is at least $\TheProportion$.
      \end{quote}
      Then $\TheGroup$ has subexponential growth.
    \end{prop}
    \begin{proof}
      Since $\AltGroup$ has finite index, there is a constant
      $\TheConstant$ such that every $\AltGroup$-coset in $\TheGroup$
      has a representative of length $\leq\TheConstant$.
      Now assume, let $\TheGrowth$ be the growth function and
      $\TheGrowthRate$ be its growth rate. Then we have, for any
      $\TheEpsilon>\Zero$ and $\TheRadius$ large enough:
      \begin{eqnarray*}
        \ToThePower{\Parentheses{
          \TheGrowthRate-\TheEpsilon
        }}{\TheRadius}
          &\leq&
          \TheGrowthOf{\TheRadius}
          \\
          &=&
          \CardOf{
            \SetOf[\TheGroupElement\in\TheGroup]{
              \TheLengthOf{\TheGroupElement} \leq \TheRadius
            }
          }
          \\
          &\leq&
          \CardOf{
            \SetOf[\AltGroupElement\in\AltGroup]{
              \TheLengthOf{\AltGroupElement} \leq \TheRadius + \TheConstant
            }
          }
          \\
          &\leq&
          \frac{\One}{\TheProportion}
          \CardOf{
            \SetOf[\AltGroupElement\in\AltGroup]{
              \TheLengthOf{\AltGroupElement} \leq \TheRadius + \TheConstant
              \text{\ and\ }
              \Sum[\TheIndex=\One][\TheNumber]{
                \TheLengthOf{
                  \TheHomomorphismOf[\TheIndex]{\AltGroupElement}
                }
              }
              \leq
              \ThePercentage\Parentheses{
                \TheRadius + \TheConstant
              }
              + \TheShiftConst
            }
          }
          \\
          &\leq&
          \frac{\One}{\TheProportion}
          \Int[{
            \substack{
              \TheLength[\One] + \cdots + \TheLength[\TheNumber] \leq
              \ThePercentage
              \Parentheses{
                \TheRadius + \TheConstant
              }
              +
              \TheShiftConst
              \\
              \TheLength[\TheIndex] \geq 0
            }
          }]{
            \TheGrowthOf{\TheLength[\One]}
            \cdots
            \TheGrowthOf{\TheLength[\TheNumber]}
            \diff \TheLength[\One]
            \cdots
            \diff \TheLength[\TheNumber]
          }
          \\
          &\leq&
          \AltPolynomialOf[\TheNumber]{\TheRadius}
          \TheGrowthRate[][
            \ThePercentage
            \Parentheses{
              \TheRadius + \TheConstant
            }
            +
            \TheShiftConst
          ]
       \end{eqnarray*}
       where $\AltPolynomial[\TheNumber]$ is a polynomial of
       degree $\TheNumber$. It follows that $\TheGrowthRate\leq\One$
       and $\TheGroup$ has subexponential growth.
    \end{proof}

  \section{Growth in Grigorchuk's First Group \boldmath$\GroupG$}
    The First Grigorchuk group is defined as
    \(
      \GroupG := \GroupPresented{
        \TheSwap,
        \TheAutom,
        \AltAutom,
        \ThrAutom
      }
    \)
    where $\TheAutom$, $\AltAutom$, and $\ThrAutom$ are binary
    tree automorphisms defined by the following system:
    \begin{eqnarray*}
      \TheAutom &=& \Pair{\TheSwap}{\AltAutom} \\
      \AltAutom &=& \Pair{\TheSwap}{\ThrAutom} \\
      \ThrAutom &=& \Pair{\TheTrivialElement}{\TheAutom}
    \end{eqnarray*}
    We start by recalling some well known facts about $\GroupG$.
    \begin{lemma}\label{KleinGruppe}
      The set $\SetOf{\TheTrivialElement,\TheAutom,\AltAutom,\ThrAutom}$
      is a subgroup of $\GroupG$ isomorphic to
      $\CyclicGroup[\Two]\crossprod\CyclicGroup[\Two]$.
    \end{lemma}
    \begin{proof}
      The system
      \begin{eqnarray*}
        \TheVar &=& \Pair{\TheTrivialElement}{\AltVar} \\
        \AltVar &=& \Pair{\TheTrivialElement}{\ThrVar} \\
        \ThrVar &=& \Pair{\TheTrivialElement}{\TheVar}
      \end{eqnarray*}
      clearly defines $\TheVar=\AltVar=\ThrVar=\TheTrivialElement$.
      On the other hand, one easily verifies that
      $\TupelOf{\TheAutom[][\Two],\AltAutom[][\Two],\ThrAutom[][\Two]}$
      and
      \(
        \TupelOf{
          \TheAutom\AltAutom\ThrAutom,
          \AltAutom\ThrAutom\TheAutom,
          \ThrAutom\TheAutom\AltAutom
        }
      \)
      are also solutions to that system. More specifically, one obtains
      the following systems from the defining equations:
      \[
        \begin{array}{rcl}
          \TheAutom[][\Two] &=& \Pair{\TheTrivialElement}{\AltAutom[][\Two]}\\
          \AltAutom[][\Two] &=& \Pair{\TheTrivialElement}{\ThrAutom[][\Two]}\\
          \ThrAutom[][\Two] &=& \Pair{\TheTrivialElement}{\TheAutom[][\Two]}
        \end{array}
        \qquad
        \qquad
        \begin{array}{rcl}
          \TheAutom\AltAutom\ThrAutom &=&
          \Pair{\TheTrivialElement}{\AltAutom\ThrAutom\TheAutom}\\
          \AltAutom\ThrAutom\TheAutom &=&
          \Pair{\TheTrivialElement}{\ThrAutom\TheAutom\AltAutom}\\
          \ThrAutom\TheAutom\AltAutom &=&
          \Pair{\TheTrivialElement}{\TheAutom\AltAutom\ThrAutom}\\
        \end{array}
      \]          
      Thus, the group generated
      by $\TheAutom$, $\AltAutom$, and $\ThrAutom$ is a quotient of
      $\CyclicGroup[\Two]\crossprod\CyclicGroup[\Two]$.
      
      Now, we see that $\TheAutom=\Pair{\TheSwap}{\AltAutom}$
      acts nontrivially on the left subtree. Thus
      $\TheAutom\neq\TheTrivialElement$. This implies
      $\ThrAutom\neq\TheTrivialElement$. Finally,
      \(
        \TheAutom=\Pair{\TheSwap}{\AltAutom}
        \neq
        \Pair{\TheTrivialElement}{\TheAutom}
        =\ThrAutom
        .
      \)
    \end{proof}
    Historically, $\GroupG$ was the first group shown to have
    intermediate growth \cite{Grigorchuk:1983}. The proof of Theorem
    \ref{GrigSubexp} included here is the result of several
    simplifications.
    \begin{theorem}[Grigorchuk]\label{GrigSubexp}
      $\GroupG$ has subexponential growth.
    \end{theorem}
    \begin{proof}
      Even elements are represented
      by words containing an even number of $\TheSwap$-letters.
      As a subgroup, $\GroupG[\One]$ is generated by the elements
      \(
        \TheAutom=\Pair{\TheSwap}{\AltAutom},
        \quad
        \AltAutom=\Pair{\TheSwap}{\ThrAutom},
        \quad
        \ThrAutom=\Pair{\TheTrivialElement}{\TheAutom},
        \quad
        \TheSwap\TheAutom\TheSwap=\Pair{\AltAutom}{\TheSwap},
        \quad
        \TheSwap\AltAutom\TheSwap=\Pair{\ThrAutom}{\TheSwap},
        \quad
        \TheSwap\TheAutom\TheSwap=\Pair{\TheAutom}{\TheTrivialElement}
        .
      \)
      As a consequence, there is an injective homomorphism
      \[
        \Pair{
          \TheHomomorphism[\Left]
        }{
          \TheHomomorphism[\Right]
        }\mapcolon
        \GroupG[\One] \monorightarrow
        \GroupG \crossprod \GroupG
        .
      \]
      Consider the length function $\TheLength$ induced by the weights
       \[
        \CoverLengthOf{\TheSwap} := \Three
        \qquad
        \CoverLengthOf{\TheAutom} := \Five
        \qquad
        \CoverLengthOf{\AltAutom} := \Four
        \qquad
        \CoverLengthOf{\ThrAutom} := \Three
        .
      \]
      We claim that
      \(
        \Pair{\TheHomomorphism[\Left]}{\TheHomomorphism[\Right]}
      \)
      and $\TheLength$
      satisfy the hypotheses of
      Proposition~\ref{BasicTool} with
      $\ThePercentage=\frac{\Seven}{\Eight}$,
      $\TheProportion=\One$, and $\TheShiftConst=\Three$.

      It follows from Lemma~\ref{KleinGruppe} that every
      minimum length word (i.e., a word that is a minimum length
      representative for the group element it represents) alternates
      between  $\TheSwap$ letters and the other letters. In other words
      a minimum length word does not contain two adjacent
      $\TheSwap$ letters
      (because they would cancel out) nor does it contain
      two adjacent non-$\TheSwap$ letters (because they would
      cancel or multiply to yield just one letter).
      Thus, in order to establish the hypotheses of
      Proposition~\ref{BasicTool}, it suffices to consider a word
      \[
        \TheWord
        =
        (\TheSwap)\TheGenerator[\One]\TheSwap\TheGenerator[\Two]
        \cdots
        \TheSwap\TheGenerator[\TheLastIndex](\TheSwap)
      \]
      with an even number of $\TheSwap$ letters (the ones at the ends
      possibly being omitted) and where
      \(
        \TheGenerator[\TheIndex]\in\SetOf{\TheAutom,\AltAutom,\ThrAutom}
        .
      \)
      Split $\TheWord$ into blocks of four letters, possibly followed by
      a single shorter block at the end. The
      homomorphism
      \(
        \Pair{\TheHomomorphism[\Left]}{\TheHomomorphism[\Right]}
      \)
      will reduce the length of each four letter block by a
      factor of at least $\frac{\Seven}{\Eight}$. For instance,
      \(
        \TheSwap\TheAutom\TheSwap\AltAutom
        \mapsto
        \Pair{\AltAutom\TheSwap}{\TheSwap\ThrAutom}
      \)
      corresponds to a reduction from length $\Number{15}$ to
      $\Number{13}$. The worst case is attained by the block
      $\TheSwap\TheAutom\TheSwap\TheAutom$ which yields a reduction from
      $\Number{16}$ to $\Number{14}$.
      Analyzing the effect of the trailing block requires a finite amount
      of checking. One obtains
      that for any element
      $\TheGroupElement\in\GroupG[\One]$,
      \[
        \TheLengthOf{
          \TheHomomorphismOf[\Left]{\TheGroupElement}
        }
        +
        \TheLengthOf{
          \TheHomomorphismOf[\Right]{\TheGroupElement}
        }
        \leq
        \frac{\Seven}{\Eight}
        \TheLengthOf{\TheGroupElement}
        +
        \Three
      \]
      as claimed, and the result follows.
    \end{proof}
    \begin{rem}
      There is nothing magic about the weights we used, and there are
      other weights that would do just as well. A particularly
      good assignment of weights was used in
      \cite{Bartholdi:2000} to obtain a good explicit
      upper bound for the growth in $\GroupG$.
    \end{rem}

  \section{Growth in \boldmath$\GroupH$}
    In this section, we will study the group
    \(
      \GroupH:=\GroupPresented{
        \TheSwap,
        \TheAutom,
        \AltAutom
      }
    \)
    where $\TheAutom$ and $\AltAutom$ are defined
    by the following equations:
    \begin{eqnarray*}
      \TheAutom &=& \Pair{\TheSwap}{\AltAutom} \\
      \AltAutom &=& \Pair{\TheTrivialElement}{\TheAutom}
    \end{eqnarray*}
    This group serves as a model to understand the slightly more involved
    group $\GroupI$ introduced in the next section. The group
    $\GroupH$ is less manageable
    than $\GroupG$ as it contains a ``self-replicating element'' of
    infinite order. For this reason, we will use the full power of
    Proposition~\ref{BasicTool} with a proportion factor
    $\TheProportion$ strictly less than $\One$.
    \begin{lemma}\label{KleinGroupH}
      The elements $\TheAutom$ and $\AltAutom$ generate a
      copy of $\CyclicGroup[\Two]\crossprod\CyclicGroup[\Two]$
      inside $\GroupH$.
    \end{lemma}
    \begin{proof}
      We have:
      \[
        \begin{array}{rcl}
          \TheAutom[][\Two] &=& \Pair{\TheTrivialElement}{\AltAutom[][\Two]}\\
          \AltAutom[][\Two] &=& \Pair{\TheTrivialElement}{\TheAutom[][\Two]}
        \end{array}
        \qquad
        \begin{array}{rcl}
          \TheAutom\AltAutom\TheAutom\AltAutom
          &=&
          \Pair{\TheTrivialElement}{\AltAutom\TheAutom\AltAutom\TheAutom}
          \\
          \AltAutom\TheAutom\AltAutom\TheAutom
          &=&
          \Pair{\TheTrivialElement}{\TheAutom\AltAutom\TheAutom\AltAutom}
        \end{array}
      \]
      Therefore:
      \[
        \TheAutom[][\Two]
        =
        \AltAutom[][\Two]
        =
        \ToThePower{
          \Parentheses{\TheAutom\AltAutom}
        }{\Two}
        =
        \TheTrivialElement
        .
      \]
      Thus $\TheAutom$ and $\AltAutom$ generate a quotient of
      $\CyclicGroup[\Two]\crossprod\CyclicGroup[\Two]$.
      Since $\TheAutom$ acts like $\TheSwap$ in the left
      subtree, we have $\TheAutom\neq\TheTrivialElement$, which
      in turn implies $\AltAutom\neq\TheTrivialElement$. Moreover,
      \(
        \TheAutom
        \neq
        \AltAutom
      \)
      by their actions on the left subtree.
      Hence $\GroupPresented{\TheAutom,\AltAutom}$ has more than
      two elements.
    \end{proof}
    \begin{prop}
      The elements $\TheSwap\TheAutom$ and $\TheSwap\AltAutom$ have
      finite order, but
      the element $\TheSwap\ThrAutom$ has infinite order.
    \end{prop}
    \begin{proof}
      First, we observe that
      \[
        \ToThePower{
          \Parentheses{\TheSwap\AltAutom}
        }{\Four}
        =
        \ToThePower{
          \Pair{\TheAutom}{\TheAutom}
        }{\Two}
        =
        \Pair{\TheTrivialElement}{\TheTrivialElement}
        =
        \TheTrivialElement
      \]
      and
      \[
        \ToThePower{
          \Parentheses{\TheSwap\TheAutom}
        }{\Two}
        =
        \Pair{\AltAutom\TheSwap}{\TheSwap\AltAutom}
        .
      \]
      Thus, $\TheSwap\AltAutom$ has order $\Four$ and
      $\TheSwap\TheAutom$ has order $\Eight$.

      On the other hand, $\ThrAutom$ replicates itself:
      \[
        \ToThePower{
          \Parentheses{
            \TheSwap\ThrAutom
          }
        }{\Two}
        =
        \TheSwap\TheAutom\AltAutom\TheSwap\TheAutom\AltAutom
        =
        \Pair{\AltAutom}{\TheSwap}\Pair{\TheAutom}{\TheTrivialElement}
        \Pair{\TheSwap}{\AltAutom}\Pair{\TheTrivialElement}{\TheAutom}
        =
        \Pair{\AltAutom\TheAutom\TheSwap}{\TheSwap\AltAutom\TheAutom}
        =
        \Pair{\ThrAutom\TheSwap}{\TheSwap\ThrAutom}
        .
      \]
      Then
      \[
        \ToThePower{
          \Parentheses{
            \TheSwap\ThrAutom
          }
        }{\Two\TheNumber}
        =
        \Pair{
          \ToThePower{
            \Parentheses{
              \ThrAutom\TheSwap
            }
          }{\TheNumber}
        }{
          \ToThePower{
            \Parentheses{
              \TheSwap\ThrAutom
            }
          }{\TheNumber}
        }
        .
      \]
      Therefore,
      \(
        \ToThePower{
          \Parentheses{
            \TheSwap\ThrAutom
          }
        }{\Two\TheNumber}
        =
        \TheTrivialElement
      \)
      implies
      \(
        \ToThePower{
          \Parentheses{
            \TheSwap\ThrAutom
          }
        }{\TheNumber}
        =
        \TheTrivialElement
        .
      \)
      So the order of $\TheSwap\ThrAutom$ is either odd or infinite.
      But any odd power of $\TheSwap\ThrAutom$ acts non-trivially on the
      tree $\TheRootedTree$ as it
      performs a swap at
      the root vertex.
    \end{proof}
    \begin{rem}
      The identity
      \[
        \ToThePower{
          \Parentheses{
            \TheSwap\ThrAutom
          }
        }{\Two\TheNumber}
        =
        \Pair{
          \ToThePower{
            \Parentheses{
              \ThrAutom\TheSwap
            }
          }{\TheNumber}
        }{
          \ToThePower{
            \Parentheses{
              \TheSwap\ThrAutom
            }
          }{\TheNumber}
        }
      \]
      rules out any hope that one could prove
      subexponential growth as easily as for
      Grigorchuk's group $\GroupG$: A block of the form
      \(
        \TheSwap\ThrAutom\TheSwap\ThrAutom
      \)
      inside an alternating word does not display any length reduction
      regardless of the choice of weights for the generators.
      Thus, there is no uniform reduction for all alternating words.
    \end{rem}
    Having located the problem, we propose a solution.

    First, we induce an appropriate length function $\TheLength$
    on $\GroupH$ via the weights
    \[
      \CoverLengthOf{\TheSwap} := \Three
      \qquad
      \CoverLengthOf{\TheAutom} := \Five
      \qquad
      \CoverLengthOf{\AltAutom} := \Four
      \qquad
      \CoverLengthOf{\ThrAutom} := \Three
      .
    \]
    Note that minimum length words alternate between the letter $\TheSwap$
    and the other letters.

    Since there is no uniform reduction, we will introduce a measure
    that allows us to distinguish between group elements that have a
    good reduction and those that do not.
    Let $\TheEpsilon$ be a small positive parameter
    to be chosen later.
    \begin{Def}
      An alternating
      word $\TheWord$ with $\TheNumber$ non-$\TheSwap$ letters
      is said to be \notion{$\TheEpsilon$-bad} if it contains at most
      $\TheEpsilon\TheNumber$ letters $\TheAutom$ and $\AltAutom$.
      An even element $\AltGroupElement\in\GroupH[\One]$ is called
      $\TheEpsilon$-bad if every minimum length representative
      of $\AltGroupElement$ is $\TheEpsilon$-bad. The element
      $\AltGroupElement$ is \notion{$\TheEpsilon$-good} if it is
      not $\TheEpsilon$-bad.
    \end{Def}
    Note that for small $\TheEpsilon$, the $\TheEpsilon$-bad \emph{words}
    will form only a slim fraction of all words of a given length. This,
    however, does not imply that the $\TheEpsilon$-bad \emph{group elements}
    are a minority among the group elements of a given length: it is
    conceivable that all the $\TheEpsilon$-bad
    words represent different group
    elements whereas the $\TheEpsilon$-good words
    represent only a handful of different
    elements.
    This is the main issue to be addressed in our proof of the following:
    \begin{prop}[Grigorchuk]
      $\GroupH$ has subexponential growth.
    \end{prop}
    \begin{proof}
      Again, we will consider the ``split into a pair'' homomorphism
      \[
        \Pair{\TheHomomorphism[\Left]}{\TheHomomorphism[\Right]}
        \mapcolon
        \GroupH[\One]
        \rightarrow
        \GroupH \crossprod \GroupH
      \]
      defined on the subgroup $\GroupH[\One]$ of even elements. Our
      goal is to verify the hypotheses of Proposition~\ref{BasicTool}
      in this setting.
    
      For the length function specified above, an alternating word
      $\TheWord$ of length $\CoverLengthOf{\TheWord}\leq\TheRadius$ contains
      at most $\Parentheses{\TheRadius+\Three}/\Six$ non-$\TheSwap$ letters.
      Let
      \(
        \TheBadWordsOf[\TheEpsilon]{\TheNumber}
      \)
      be the number of $\TheEpsilon$-bad alternating words that have
      $\TheNumber$ non-$\TheSwap$ letters and represent an even element
      of $\GroupH$.
      Then
      \[
        \TheBadWordsOf[\TheEpsilon]{\TheNumber}
        \leq
        \Two \Choose{\TheNumber}{\FloorOf{\TheEpsilon\TheNumber}}
        \ToThePower{\Three}{\FloorOf{\TheEpsilon\TheNumber}}
        ,
      \]
      where he leading $\Two$ accounts for the choice of starting $\TheWord$
      with $\TheSwap$ or not, the binomial coefficient counts the number of
      ways to allocate
      $\FloorOf{\TheEpsilon\TheNumber}$ positions where
      $\TheAutom$ or $\AltAutom$ can be placed, and the factor
      $\ToThePower{\Three}{\FloorOf{\TheEpsilon\TheNumber}}$ counts the
      ways of using the allocated positions: note that
      the letter $\ThrAutom$ can be used in any of these locations.

      We estimate the number
      \(
        \TheBadElementsOf[\TheEpsilon]{\TheRadius}
      \)
      of $\TheEpsilon$-bad group elements
      of length at most $\TheRadius$ by
      \[
        \TheBadElementsOf[\TheEpsilon]{\TheRadius}
        \leq
        \Sum[\TheNumber=\Zero][\Parentheses{\TheRadius+\Three}/\Six]{
          \TheBadWordsOf[\TheEpsilon]{\TheNumber}
        }
        \leq
        \frac{\Parentheses{\TheRadius+\Three}}{\Three}
        \Choose{\Parentheses{\TheRadius+\Three}/\Six}{\FloorOf{\TheEpsilon\Parentheses{\TheRadius+\Three}/\Six}}
        \Three[][\FloorOf{\TheEpsilon\Parentheses{\TheRadius+\Three}/\Six}]
      \]
      
      Now, suppose $\GroupH$ has exponential growth. Then the index~$\Two$
      subgroup $\GroupH[\One]$ has exponential growth with respect to
      the restricted length
      \(
        \TheLength\RestrictionTo{\GroupH[\One]}
        .
      \)
      Let $\TheGrowthRate>\One$ be the growth rate of $\GroupH[\One]$.
      Using Stirling's formula, we can choose $\TheEpsilon>\Zero$
      small enough
      so that $\TheGrowthRate[][\TheRadius]$ grows faster than
      $\TheBadElementsOf[\TheEpsilon]{\TheRadius}$. In particular,
      for $\TheRadius$ big enough, at least half of the elements in
      \(
        \SetOf[{
          \AltGroupElement\in\GroupH[\One]
        }]{
          \TheLengthOf{\AltGroupElement} \leq \TheRadius
        }
      \)
      are $\TheEpsilon$-good.

      We finish the proof by observing that an $\TheEpsilon$-good
      element reduces at least by a factor of
      \[
        \ThePercentage :=
        \frac{
          \Four\TheEpsilon + \Parentheses{\Two-\TheEpsilon}\Three
        }{
          \Five\TheEpsilon + \Parentheses{\Two-\TheEpsilon}\Three
        }
        <
        \One
      \]
      when written in the pair notation.
      The worst case is attained when among the
      non-$\TheSwap$ letters, we have the largest possible number of
      $\ThrAutom$ letters and all the remaining are $\TheAutom$.
      We have thus verified the hypotheses of
      Proposition~\ref{BasicTool} with $\ThePercentage < \One$ and
      $\TheProportion = \One/\Two$.
    \end{proof}
    \begin{rem}
      Note that in this proof, we had to fix the assumed
      growth rate $\TheGrowthRate>\One$ in order to find an
      $\TheEpsilon>\Zero$ small enough
      to verify the hypotheses of Proposition~\ref{BasicTool}. This
      method will not allow us to deduce a bound of the form
      \[
        \TheGrowthOf{\TheRadius}
        \leq
        \ExpOf{\AltPercentage\TheRadius}
      \]
      with $\AltPercentage<\One$. Such a bound can be deduced in the
      case of Grigorchuk's group since the reduction factor
      $\frac{\Seven}{\Eight}$ holds independently of the assumed growth rate.
      Details can be found in \cite{Bartholdi:2000}.
    \end{rem}

  \section{Growth in \boldmath$\GroupI$}
    The critical orbit of the polynomial $\TheComplex[][\Two] + \ImUnit$
    is
    \[
      \Zero \mapsto
      \ImUnit \mapsto
      \Parentheses{-\One+\ImUnit} \mapsto
      -\ImUnit \mapsto
      \Parentheses{-\One+\ImUnit}
    \]
    i.e., $\TheComplex[][\Two] + \ImUnit$
    is postcritically finite. As a consequence, the iterated
    monodromy group
    \(
      \GroupI := \ImgOf{\TheComplex[][\Two] + \ImUnit}
    \)
    is generated by three elements descended from loops around
    \(
      \ImUnit,
    \)
    \(
      \Parentheses{-\One+\ImUnit}
      ,
    \)
    and
    \(
      -\ImUnit
      .
    \)
    More precisely, according to
    \cite{Bartholdi.Grigorchuk.Nekrashevych:2002},
    $\GroupI=\GroupPresented{\TheSwap,\TheAutom,\AltAutom}$
    where $\TheAutom$,
    and $\AltAutom$ are defined by the following equations:
    \begin{eqnarray*}
      \TheAutom &=& \Pair{\TheSwap}{\AltAutom} \\
      \AltAutom &=& \Pair{\TheAutom}{\TheTrivialElement}
    \end{eqnarray*}
    The generators $\TheAutom$
    and $\AltAutom$ correspond to loops around the periodic points
    $\Parentheses{-\One+\ImUnit}$
    and $-\ImUnit$, while $\TheSwap$ represents the swap
    on $\TheRootedTree$, which is induced by the loop around the preperiodic
    point $\ImUnit$.

    The equations defining $\GroupI$ are very similar to those for
    $\GroupH$ and yet, a new complication appears: now we are facing
    several non-reducing elements.
    \begin{lemma}\label{DihedralGroupI}
      Any two of the generators span a finite dihedral group
      inside $\GroupI$:
      \begin{eqnarray*}
        \DihedralGroupOf[\Four]{\TheSwap,\AltAutom}
        &:=&
        \GroupPresented{
          \TheSwap,\AltAutom
        }
        =
        \GroupPresented[
          \TheSwap,\AltAutom
        ]{
          \TheSwap[][\Two]
          =
          \AltAutom[][\Two]
          =
          \ToThePower{
            \Parentheses{
              \TheSwap
              \AltAutom
            }
          }{\Four}
          =
          \TheTrivialElement
        }
        \\
        \DihedralGroupOf[\Eight]{\TheSwap,\TheAutom}
        &:=&
        \GroupPresented{
          \TheSwap,\TheAutom
        }
        =
        \GroupPresented[
          \TheSwap,\TheAutom
        ]{
          \TheSwap[][\Two]
          =
          \TheAutom[][\Two]
          =
          \ToThePower{
            \Parentheses{
              \TheSwap
              \TheAutom
            }
          }{\Eight}
          =
          \TheTrivialElement
        }
        \\
        \DihedralGroupOf[\Eight]{\TheAutom,\AltAutom}
        &:=&
        \GroupPresented{
          \TheAutom,\AltAutom
        }
        =
        \GroupPresented[
          \TheAutom,\AltAutom
        ]{
          \TheAutom[][\Two]
          =
          \AltAutom[][\Two]
          =
          \ToThePower{
            \Parentheses{
              \TheAutom
              \AltAutom
            }
          }{\Eight}
          =
          \TheTrivialElement
        }
      \end{eqnarray*}
    \end{lemma}
    \begin{proof}
      Note that $\TheAutom$ is not trivial since it acts like $\TheSwap$
      in the left subtree. Hence
      \(
        \AltAutom
        =
        \Pair{\TheAutom}{\TheTrivialElement}
      \)
      is not trivial either.
      Since $\TheAutom[][\Two]=\Pair{\TheTrivialElement}{\AltAutom[][\Two]}$
      and $\AltAutom[][\Two]=\Pair{\TheAutom[][\Two]}{\TheTrivialElement}$,
      we infer that $\TheAutom$ and $\AltAutom$ have order $\Two$.

      From
      \(
        \TheSwap\AltAutom\TheSwap\AltAutom
        =
        \Pair{\TheAutom}{\TheAutom}
      \)
      We infer that $\TheSwap\AltAutom$ (and its inverse $\AltAutom\TheSwap$)
      have order $\Four$. Hence
      \(
        \TheSwap\TheAutom\TheSwap\TheAutom
        =
        \Pair{\AltAutom\TheSwap}{\TheSwap\AltAutom}
      \)
      implies that $\TheSwap\TheAutom$ has order $\Eight$. Thus
      $\TheAutom\AltAutom=
      \Pair{\TheSwap\TheAutom}{\AltAutom}$ has order $\Eight$, too.
    \end{proof}

    We want to ensure that every group element is represented by
    a  minimum length word that alternates between the letter $\TheSwap$
    and other generators.
    Using Lemma~\ref{DihedralGroupI}, we achieve this goal by
    adding the non-trivial elements of
    $\DihedralGroupOf[\Eight]{\TheAutom,\AltAutom}$ to
    the designated generating set for $\GroupI$. For reference, we
    list the extended set of generators:
    \[
      \begin{array}{c|c|c}
        \text{name} & \text{element} & \text{pair} \\
        \TheAutom[\One] & \TheAutom & \Pair{\TheSwap}{\AltAutom} \\ 
        \TheAutom[\Two] & \TheAutom\AltAutom & \Pair{\TheSwap\TheAutom}{\AltAutom} \\ 
        \TheAutom[\Three] &\TheAutom\AltAutom\TheAutom & \Pair{\TheSwap\TheAutom\TheSwap}{\TheTrivialElement} \\ 
        \TheAutom[\Four] &\TheAutom\AltAutom\TheAutom\AltAutom & \Pair{\TheSwap\TheAutom\TheSwap\TheAutom}{\TheTrivialElement} \\ 
        \TheAutom[\Five] &\TheAutom\AltAutom\TheAutom\AltAutom\TheAutom & \Pair{\TheSwap\TheAutom\TheSwap\TheAutom\TheSwap}{\AltAutom} \\ 
        \TheAutom[\Six] &\TheAutom\AltAutom\TheAutom\AltAutom\TheAutom\AltAutom & \Pair{\TheSwap\TheAutom\TheSwap\TheAutom\TheSwap\TheAutom}{\AltAutom} \\ 
        \TheAutom[\Seven] &\TheAutom\AltAutom\TheAutom\AltAutom\TheAutom\AltAutom\TheAutom & \Pair{\TheSwap\TheAutom\TheSwap\TheAutom\TheSwap\TheAutom\TheSwap}{\TheTrivialElement} \\ 
        \TheAutom[\Eight] &\TheAutom\AltAutom\TheAutom\AltAutom\TheAutom\AltAutom\TheAutom\AltAutom & \Pair{\TheSwap\TheAutom\TheSwap\TheAutom\TheSwap\TheAutom\TheSwap\TheAutom}{\TheTrivialElement} \\ 
      \end{array}
      \qquad
      \begin{array}{c|c|c}
        \text{name} & \text{element} & \text{pair} \\ 
        \AltAutom[\One] & \AltAutom & \Pair{\TheAutom}{\TheTrivialElement} \\ 
        \AltAutom[\Two] & \AltAutom\TheAutom & \Pair{\TheAutom\TheSwap}{\AltAutom} \\ 
        \AltAutom[\Three] &\AltAutom\TheAutom\AltAutom & \Pair{\TheAutom\TheSwap\TheAutom}{\AltAutom} \\ 
        \AltAutom[\Four] &\AltAutom\TheAutom\AltAutom\TheAutom & \Pair{\TheAutom\TheSwap\TheAutom\TheSwap}{\TheTrivialElement} \\ 
        \AltAutom[\Five] &\AltAutom\TheAutom\AltAutom\TheAutom\AltAutom & \Pair{\TheAutom\TheSwap\TheAutom\TheSwap\TheAutom}{\TheTrivialElement} \\ 
        \AltAutom[\Six] &\AltAutom\TheAutom\AltAutom\TheAutom\AltAutom\TheAutom & \Pair{\TheAutom\TheSwap\TheAutom\TheSwap\TheAutom\TheSwap}{\AltAutom} \\ 
        \AltAutom[\Seven] &\AltAutom\TheAutom\AltAutom\TheAutom\AltAutom\TheAutom\AltAutom & \Pair{\TheAutom\TheSwap\TheAutom\TheSwap\TheAutom\TheSwap\TheAutom}{\AltAutom} \\ 
        \AltAutom[\Eight] &\AltAutom\TheAutom\AltAutom\TheAutom\AltAutom\TheAutom\AltAutom\TheAutom & \Pair{\TheAutom\TheSwap\TheAutom\TheSwap\TheAutom\TheSwap\TheAutom\TheSwap}{\TheTrivialElement} \\ 
      \end{array}
    \]
    Note that $\TheAutom[\Eight]=\AltAutom[\Eight]$.

    \begin{rem}
      The issue of choosing weights for the generators is delicate
      since we will insist that words do not increase in length when
      split into pair notation.
      We will choose weights only for $\TheSwap$,
      $\TheAutom$, and $\AltAutom$; the weights for the
      redundant generators will simply be their lengths relative to
      the smaller generating set.

      Note that the equation
      \(
        \TheSwap\TheAutom\TheSwap\TheAutom
        =
        \Pair{
          \AltAutom\TheSwap
        }{
          \TheSwap\AltAutom
        }
      \)
      imposes the restriction
      \(
        \CoverLengthOf{\TheAutom} \geq \CoverLengthOf{\AltAutom}
      \)
      while
      \(
        \TheSwap\AltAutom\TheAutom\AltAutom
        \TheSwap\AltAutom\TheAutom\AltAutom
        =
        \Pair{
          \AltAutom\TheAutom\TheSwap\TheAutom
        }{
          \TheAutom\TheSwap\TheAutom\AltAutom
        }
      \)
      requires
      \(
        \CoverLengthOf{\AltAutom} \geq \CoverLengthOf{\TheAutom}
        .
      \)
      We will settle upon the weights
      \[
        \CoverLengthOf{\TheSwap} = \Three,\qquad
        \CoverLengthOf{\TheAutom} = \Four,\qquad
        \CoverLengthOf{\AltAutom} = \Four
      \]
      and denote the induced length function by $\TheLength$.
    \end{rem}

    As before, we will use the homomorphism
    \[
      \Pair{\TheHomomorphism[\Left]}{\TheHomomorphism[\Right]}
      \mapcolon
      \GroupI[\One]
      \rightarrow
      \GroupI \crossprod \GroupI
    \]
    that splits even elements. This time, however, we will have to
    split elements up to three times. Note that every generator
    \(
      \TheGenerator \in \DihedralGroupOf[\Eight]{\TheAutom,\AltAutom}
    \)
    satisfies
    \(
      \TheLengthOf{\TheHomomorphismOf[\Left]{\TheGenerator}}
      +
      \TheLengthOf{\TheHomomorphismOf[\Right]{\TheGenerator}}
      \leq
      \TheLengthOf{\TheSwap}
      +
      \TheLengthOf{\TheGenerator}
      .
    \)
    This means that any length reduction secured in the process of
    splitting a word into a pair will persist in subsequent splittings.
    We call a generator
    \(
      \TheGenerator\in\DihedralGroupOf[\Eight]{\TheAutom,\AltAutom}
      \setminus\SetOf{\TheTrivialElement}
    \)
    \notion{good by nature}
    if the inequality above is strict:
    \[
      \TheLengthOf{\TheHomomorphismOf[\Left]{\TheGenerator}}
      +
      \TheLengthOf{\TheHomomorphismOf[\Right]{\TheGenerator}}
      <
      \TheLengthOf{\TheSwap} + \TheLengthOf{\TheGenerator}
      .
    \]
    Otherwise, $\TheGenerator$ is \notion{bad}.
    \begin{observation}
      The only bad letters are $\TheAutom$, $\TheAutom[\Two]$,
      $\AltAutom[\Two]$, and $\AltAutom[\Three]$.
      Any generator $\TheGenerator$ that is good by nature satisfies
      \(
        \TheLengthOf{\TheHomomorphismOf[\Left]{\TheGenerator}}
        +
        \TheLengthOf{\TheHomomorphismOf[\Right]{\TheGenerator}}
        \leq
        \frac{\Number{29}}{\Number{31}}
        \Parentheses{
          \TheLengthOf{\TheSwap} + \TheLengthOf{\TheGenerator}
        }
        .
      \)
      The worst case is attained by $\AltAutom[\Seven]$
      since
      \(
        \TheLengthOf{\TheHomomorphismOf[\Left]{\AltAutom[\Seven]}}
        +
        \TheLengthOf{\TheHomomorphismOf[\Right]{\AltAutom[\Seven]}}
        = \Number{29}
      \)
      while
      \(
        \TheLengthOf{\TheSwap}
        +
        \TheLengthOf{\AltAutom[\Seven]}
        = \Number{31}
        .
      \)
      As a consequence, if the even group element
      \(
        \TheGroupElement
        =
        \Pair
          {\TheHomomorphismOf[\Left]{\TheGroupElement}}
          {\TheHomomorphismOf[\Right]{\TheGroupElement}}
      \)
      can be represented by an alternating minimum length word $\TheWord$
      such that all non-$\TheSwap$ letters in $\TheWord$ are
      good by nature, then
      \[
        \TheLengthOf{\TheHomomorphismOf[\Left]{\TheGroupElement}}
        +
        \TheLengthOf{\TheHomomorphismOf[\Right]{\TheGroupElement}}
        \leq
        \frac{\Number{29}}{\Number{31}}
        \TheLengthOf{\TheWord}
        .
      \]
    \end{observation}
    The presence of several bad generators poses a problem for the counting
    method employed in the previous section: The count of $\TheEpsilon$-bad
    words of a given size had a
    factor $\Three[][\FloorOf{\TheEpsilon\TheNumber}]$
    accounting for the ways of writing letters (good or bad) at the
    selected positions. All other positions were filled with the bad
    letter $\ThrAutom$. Now that there are four bad letters, we
    would have to include an additional factor of
    $\Four[][\TheNumber - \FloorOf{\TheEpsilon\TheNumber}]$.
    This clearly grows too fast.

    To overcome this problem, we will follow the slogan ``Some
    bad letters are \notion{good by
    position}''. The crucial insight is that some combinations of
    bad letters produce good letters after splitting.

      Our method is based on the following two basic patterns:
      Let each occurrence of $\Slot$ represent an arbitrary bad
      letter, and let $\TheWord$ be a minimum length word, not
      necessarily alternating. Note that
      \(
        \TheHomomorphismOf[\Right]{\Slot} = \AltAutom
        .
      \)
      \begin{propositionsystem}{P}
        \item\label{PatternOne}
          Every substring $\TheAutom\TheSwap\Slot\TheSwap\TheAutom$
          in $\TheWord$ gives rise to the string $\TheSwap\AltAutom\TheSwap$
          in either $\TheHomomorphismOf[\Left]{\TheWord}$ or
          $\TheHomomorphismOf[\Right]{\TheWord}$ (depending on the
          position of
          $\TheAutom\TheSwap\Slot\TheSwap\TheAutom$ within $\TheWord$).
        \item\label{PatternTwo}
          Every substring $\AltAutom\TheSwap\Slot\TheSwap\AltAutom$
          in $\TheWord$ gives rise to the
          string $\TheAutom\AltAutom\TheAutom$
          in either $\TheHomomorphismOf[\Left]{\TheWord}$ or
          $\TheHomomorphismOf[\Right]{\TheWord}$. Possibly joining with
          adjacent letters, $\TheAutom\AltAutom\TheAutom$ evaluates to a
          (good) generator $\TheAutom[\TheIndex]$
          with $\TheIndex \geq \Three$ or $\AltAutom[\AltIndex]$ with
          $\AltIndex \geq \Four$.
      \end{propositionsystem}
      We now turn the basic patterns into a sufficient list of
      length reducing patterns.
      Let $\TheWord$ be an even word. The following claims
      are easily verified by using the patterns
      \pref{PatternOne} and \pref{PatternTwo}.
      \begin{propositionsystem}{B}
        \item\label{PatternA}
          If
          \(
            \TheWord \in \SetOf{
              \TheAutom\TheSwap\Slot\TheSwap\TheAutom,
              \TheAutom\TheSwap\Slot\TheSwap\TheAutom[\Two],
              \AltAutom[\Two]\TheSwap\Slot\TheSwap\TheAutom,
              \AltAutom[\Two]\TheSwap\Slot\TheSwap\TheAutom[\Two]
            }
          \)
          then
          \(
            \TheHomomorphismOf[\Left]{\TheWord} =
            \TheVariable\TheSwap\AltAutom\TheSwap\AltVariable
          \)
          with
          \(
            \TheVariable,\AltVariable\in
            \DihedralGroupOf[\Eight]{\TheAutom,\AltAutom}
            .
          \)
          Thus, $\TheHomomorphismOf[\Left]{\TheWord}$ contains a
          good generator.
        \item\label{PatternB}
          If 
          \(
            \TheWord \in \SetOf{
              \TheAutom[\Two]\TheSwap\Slot\TheSwap\AltAutom[\Two],
              \TheAutom[\Two]\TheSwap\Slot\TheSwap\AltAutom[\Three],
              \AltAutom[\Three]\TheSwap\Slot\TheSwap\AltAutom[\Two],
              \AltAutom[\Three]\TheSwap\Slot\TheSwap\AltAutom[\Three]
            }
          \)
          then
          \(
            \TheHomomorphismOf[\Left]{\TheWord}
            =
            \TheVariable\TheSwap\TheAutom[\Three]\TheSwap\AltVariable
          \)
          with
          \(
            \TheVariable,\AltVariable\in
            \DihedralGroupOf[\Eight]{\TheAutom,\AltAutom}
            .
          \)
          Thus, $\TheHomomorphismOf[\Left]{\TheWord}$ contains a
          good generator.
        \item\label{PatternC}
          If
          \(
            \TheWord \in \SetOf{
              \TheSwap\Slot\TheSwap\AltAutom[\Two]\TheSwap\Slot\TheSwap\AltAutom[\Three],
              \AltAutom[\Two]\TheSwap\Slot\TheSwap\AltAutom[\Three]\TheSwap\Slot\TheSwap
            }
          \)
          then
          \(
            \TheHomomorphismOf[\Left]{\TheWord}
            \in
            \SetOf{
              \AltAutom[\Two]\TheSwap\AltAutom[\Two]\TheSwap\TheAutom
              ,
              \TheAutom\TheSwap\AltAutom[\Two]\TheSwap\TheAutom[\Two]
            }
            .
          \)
          Thus, after splitting, this realizes a pattern from
          \pref{PatternA}.
        \item\label{PatternD}
          If
          \(
            \TheWord \in
            \SetOf{
              \Slot\TheSwap\TheAutom[\Two]\TheSwap\Slot\TheSwap\TheAutom\TheSwap\Slot,
              \Slot\TheSwap\TheAutom\TheSwap\Slot\TheSwap\AltAutom[\Two]\TheSwap\Slot
            }
          \)
          then
          \(
            \TheHomomorphismOf[\Right]{\TheWord}
            \in
            \SetOf{
              \AltAutom\TheSwap\TheAutom[\Two]\TheSwap\AltAutom,
              \AltAutom\TheSwap\AltAutom[\Two]\TheSwap\AltAutom
            },
          \)
          which realizes a pattern from \pref{PatternB}.
      \end{propositionsystem}
    \begin{rem}
      Let us illustrate how \pref{PatternOne} and \pref{PatternTwo} serve
      as the design principles for the
      list (\ref{PatternA}--\ref{PatternD}). The patterns in
      \pref{PatternA} are directly modeled upon \pref{PatternOne}; for
      instance:
      \(
        \TheAutom\TheSwap\Slot\TheSwap\TheAutom[\Two]
        =
        \TheAutom\TheSwap\Slot\TheSwap\TheAutom\AltAutom
      \)
      features $\TheAutom\TheSwap\Slot\TheSwap\TheAutom$ as a subword.
      Similarly, \pref{PatternB} relies on \pref{PatternTwo}.

      The patterns in \pref{PatternC} descend from \pref{PatternOne}:
      The common substring
      \(
        \AltAutom[\Two]\TheSwap\Slot\TheSwap\AltAutom[\Three]
      \)
      splits as follows
      \(
        \TheHomomorphismOf[\Left]{
          \AltAutom[\Two]\TheSwap\Slot\TheSwap\AltAutom[\Three]
        }
        =
        \TheAutom\TheSwap\Slot\TheSwap\TheAutom
        .
      \)
      Thus we recover \pref{PatternOne}.
      A similar computation explains how \pref{PatternD} is
      modeled after \pref{PatternTwo}.
    \end{rem}
      
    If a block of consecutive letters in a word $\TheWord$
    matches one of the patterns \pref{PatternA}--\pref{PatternD},
    it is called a 
    \notion{good block}. A letter within $\TheWord$ that is bad
    by nature is called \notion{good by position} if it belongs to
    a good block. Blocks conforming to \pref{PatternA} or
    \pref{PatternB} yield an occurrence of a good generator in
    either $\TheHomomorphismOf[\Left]{\TheWord}$ or
    $\TheHomomorphismOf[\Right]{\TheWord}$, depending on the
    position of the block inside $\TheWord$. In either case, we will
    find a reduction in length after splitting $\TheWord$ a second time
    into
    \(
      \TupelOf{
        \TheHomomorphismOf[\Left\Left]{\TheWord},
        \TheHomomorphismOf[\Left\Right]{\TheWord},
        \TheHomomorphismOf[\Right\Left]{\TheWord},
        \TheHomomorphismOf[\Right\Right]{\TheWord}
      }
      .
    \)
    Similarly, blocks conforming to \pref{PatternC} or
    \pref{PatternD} display
    a length reduction at the third splitting
    \(
      \TupelOf{
        \TheHomomorphismOf[\Left\Left\Left]{\TheWord},\cdots,
        \TheHomomorphismOf[\Right\Right\Right]{\TheWord}
      }
      .
    \)

    \begin{lemma}\label{BadStrings}
      For any number $\TheIndex$ let $\BadStrings[\TheIndex]$
      be the set of all alternating words $\TheWord$ satisfying the following
      requirements:
      \begin{itemize}
        \item
          The number of non-$\TheSwap$ letters in $\TheWord$ is
          $\TheIndex$.
        \item
          All non-$\TheSwap$ letters in $\TheWord$ are bad.
        \item
          The word $\TheWord$ does not contain good blocks.
      \end{itemize}
      Then, there is a global bound $\TheBound$ such that
      \(
        \CardOf{\BadStrings[\TheIndex]} \leq \TheBound
      \)
      for all $\TheIndex$.
    \end{lemma}
    \begin{proof}
      Consider a block $\TheVariable\TheSwap\Slot\TheSwap\AltVariable$
      in $\TheWord$. Out of all sixteen possible combinations for
      \(
        \TheVariable,\AltVariable
        \in
        \SetOf{
          \TheAutom,\TheAutom[\Two],
          \AltAutom[\Two],\AltAutom[\Three]
        },
      \)
      the twelve listed in (\ref{PatternA}--\ref{PatternD}) give rise
      to good blocks, and consequently are not featured in $\TheWord$.
      Thus, there are only four possible patterns left that can occur
      in $\TheWord$, namely:
      \(
        \TheAutom\TheSwap\Slot\TheSwap\AltAutom[\Three],
        \TheAutom[\Two]\TheSwap\Slot\TheSwap\TheAutom[\Two],
        \AltAutom[\Two]\TheSwap\Slot\TheSwap\AltAutom[\Two],
        \AltAutom[\Three]\TheSwap\Slot\TheSwap\TheAutom
        .
      \)
      Thus, we find that besides minor disturbances near the ends of
      $\TheWord$, any two positions congruent modulo $\Eight$ will
      feature the same letter.
    \end{proof}

    \begin{theorem}
      $\GroupI$ has subexponential growth.
    \end{theorem}
    \begin{proof}
      We will verify the hypotheses of Proposition~\ref{BasicTool} for the
      homomorphism
      \begin{eqnarray*}
        \GroupI[\Three]
        & \rightarrow &
        \overbrace{\GroupI \crossprod \cdots \crossprod \GroupI}^{\Eight\text{ times}} \\
        \TheGroupElement &\mapsto &
          \TupelOf{
            \TheHomomorphismOf[\Left\Left\Left]{\TheGroupElement}
            ,\cdots,
            \TheHomomorphismOf[\Right\Right\Right]{\TheGroupElement}
          }
      \end{eqnarray*}
      
      For small $\TheEpsilon>\Zero$, let
      $\TheBadWordsOf[\TheEpsilon]{\TheNumber}$
      count the number of alternating words
      representing elements in $\GroupI[\Three]$ and such that at most
      $\TheEpsilon\TheNumber$ of its $\TheNumber$ non-$\TheSwap$ letters are
      good, either by nature or by position. We have
      \[
        \TheBadWordsOf[\TheEpsilon]{\TheNumber}
        \leq
        \Two
        \Choose{\TheNumber}{\FloorOf{\TheEpsilon\TheNumber}}
        \ToThePower{\Number{15}}{\FloorOf{\TheEpsilon\TheNumber}}
        \ToThePower{\TheBound}{\One+\FloorOf{\TheEpsilon\TheNumber}}
      \]
      where the factors arise as follows:
      \begin{itemize}
        \item
          The leading $\Two$ accounts for the choice of starting
          with the letter $\TheSwap$ or not.
        \item
          The binomial coefficient selects the positions for possible
          good letters.
        \item
          The power of $\Number{15}$
          counts the ways of placing non-$\TheSwap$ letters at the selected
          positions.
        \item
          The power of $\TheBound$ accounts for the fact that the selected
          positions break the word into
          $\One+\FloorOf{\TheEpsilon\TheNumber}$
          complementary components. These subwords do not contain good
          letters, i.e., they belong to the sets $\BadStrings[\TheIndex]$
          of Lemma~\ref{BadStrings}.
      \end{itemize}

      Since the longest generator
      $\TheAutom[\Eight]$ has length $\Number{32}$, a word of length
      $\TheRadius$ has at most
      \(
        \Parentheses{\TheRadius+\Number{32}}/\Number{35}
      \)
      non-$\TheSwap$ letters. It follows that the number
      $\TheBadElementsOf[\TheEpsilon]{\TheRadius}$ of $\TheEpsilon$-bad
      elements in
      \(
        \SetOf[{\TheGroupElement\in\GroupI[\Three]}]{
          \TheLengthOf{\TheGroupElement} \leq \TheRadius
        }
      \)
      is bounded from above by
      \[
        \TheBadElementsOf[\TheEpsilon]{\TheRadius}
        \leq
        \Sum[\TheNumber=0][\Parentheses{\TheRadius+\Number{32}}/\Number{35}]{
          \TheBadWordsOf[\TheEpsilon]{\TheNumber}
        }
        .
      \]
      Assume, by contradiction, that $\GroupI[\Three]$ has
      exponential growth, with respect to the restricted
      length function $\TheLength\RestrictionTo{\GroupI[\Three]}$, and
      let $\TheGrowthRate>\One$ be the corresponding growth rate.
      Choose $\TheEpsilon>\Zero$ so
      that $\TheGrowthRate[][\TheRadius]$ grows faster than the number
      of $\TheEpsilon$-bad
      elements $\TheBadElementsOf[\TheEpsilon]{\TheRadius}$.

      The proof ends with the observation that good letters
      (by nature or by position) yield a definite reduction after
      splitting three times. Thus there is a number $\ThePercentage$
      depending only on $\TheEpsilon$ satisfying the hypotheses
      of Proposition~\ref{BasicTool}.
    \end{proof}

  \section*{References}
    \InputIfFileExists{bibliography.tex}{}{}
    
\end{document}